\documentclass[11pt]{amsart}

\usepackage{amsmath}
\usepackage{amssymb}
\usepackage{latexsym}
{\theoremstyle{plain}%
  \newtheorem{theorem}{Theorem}[section]
  \newtheorem{corollary}[theorem]{Corollary}
  \newtheorem{proposition}[theorem]{Proposition}
  \newtheorem{lemma}[theorem]{Lemma}

{\theoremstyle{remark}

\newtheorem{remark}[theorem]{Remark}
}

{\theoremstyle{definition}
\newtheorem{definition}[theorem]{Definition}
} 

\begin{document}

\newcommand{\pntwo}{\mathbb{P}^{n_2} \times \cdots \times \mathbb{P}^{n_k}}
\newcommand{\bj}{\underline{j}}
\newcommand{\bi}{\underline{i}}
\newcommand{\dep}{\operatorname{depth}}
\newcommand{\hit}{\operatorname{ht}}
\newcommand{\opdd}{\overline{P}_{d_1,d_2}}
\newcommand{\opdk}{\overline{P}_{d_1,\ldots,d_k}}
\newcommand{\op}{\overline{P}}
\newcommand{\xab}{X^{\underline{a}}Y^{\underline{b}}}
\newcommand{\xka}{X_1^{\underline{\alpha}_1}\cdots X_k^{\underline{\alpha}_k}}
\newcommand{\bb}{{\bf b}}
\newcommand{\ab}{(\underline{a},\underline{b})}
\newcommand{\alphak}{(\underline{\alpha}_1,\ldots,\underline{\alpha}_k)}
\newcommand{\abo}{(\overline{\underline{a},\underline{b}})}
\newcommand{\alphako}{\overline{(\underline{\alpha}_1,\ldots,\underline{\alpha}_k)}}
\newcommand{\opn}{\mathcal{O}_{\pr^n}}
\newcommand{\opthree}{\mathcal{O}_{\pr^3}}
\newcommand{\Iz}{I_{\Z}}
\newcommand{\Ixp}{I_{\xp}}
\newcommand{\Z}{\mathbb{Z}}
\newcommand{\xpi}{\X_{P_i}}
\newcommand{\xq}{\X_{Q_1}}
\newcommand{\xqi}{\X_{Q_i}}
\newcommand{\xp}{\X_{P_1}}
\newcommand{\lp}{L_{P_1}}
\newcommand{\ax}{\alpha_{\X}}
\newcommand{\bx}{\beta_{\X}}
\newcommand{\kxo}{k[x_0,\ldots,x_n]}
\newcommand{\kx}{k[x_1,\ldots,x_n]}
\newcommand{\popo}{\mathbb{P}^1 \times \mathbb{P}^1}
\newcommand{\pr}{\mathbb{P}}
\newcommand{\pn}{\mathbb{P}^n}
\newcommand{\pnpm}{\mathbb{P}^n \times \mathbb{P}^m}
\newcommand{\pnk}{\mathbb{P}^{n_1} \times \cdots \times \mathbb{P}^{n_k}}
\newcommand{\X}{\mathbb{X}}
\newcommand{\Y}{\mathbb{Y}}
\newcommand{\N}{\mathbb{N}}
\newcommand{\M}{\mathbb{M}}
\newcommand{\Q}{\mathbb{Q}}
\newcommand{\Ix}{I_{\X}}
\newcommand{\pix}{\pi_1(\X)}
\newcommand{\pixt}{\pi_2(\X)}
\newcommand{\pipi}{\pi_1^{-1}(P_i)}
\newcommand{\piqi}{\pi_2^{-1}(Q_i)}
\newcommand{\qpi}{Q_{P_i}}
\newcommand{\pqi}{P_{Q_i}}
\newcommand{\pitk}{\pi_{2,\ldots,k}}
\newcommand{\kdim}{\operatorname{K-}\dim}


\title{The Hilbert Functions of 
ACM Sets of Points in $\pnk$} 
\thanks{Updated: Feb. 25, 2002}

\author{Adam Van Tuyl}
\address{Department of Mathematical Sciences \\ 
Lakehead University \\ 
Thunder Bay, ON P7B 5E1, Canada}
\email{avantuyl@sleet.lakeheadu.ca}

\begin{abstract}  If $\X$ is a set of points in $\pnk$,  then the
associated coordinate ring $R/\Ix$ is an $\N^k$-graded ring.
The Hilbert function of $\X$, defined by
$H_{\X}(\bi) := \dim_{\bf k} (R/\Ix)_{\bi}$
for all $\bi \in \N^k$, is studied.   Since the ring $R/\Ix$ may or may not 
be Cohen-Macaulay, we consider only those $\X$ that
are ACM.   Generalizing
the case of $k = 1$ to all $k$, we show
that a function is the Hilbert function of
an ACM set of points if and only if its first
difference function is the Hilbert function of a multi-graded
artinian quotient.  We also give a new characterization of ACM
sets of points in $\popo$, and show how the graded Betti numbers
(and hence, Hilbert function) of ACM sets of points in this space can be
obtained solely through combinatorial means.
\end{abstract}

\maketitle

\section{Introduction}
Let $R
= {\bf k}[x_{1,0},\ldots,x_{1,n_1},\ldots,x_{k,0},\ldots,x_{k,n_k}]$
with $\deg x_{i,j} = e_i$ be the $\N^k$-graded coordinate
ring associated to $\pnk$.  A point $P = \mathcal{P}_1 \times
\cdots \times \mathcal{P}_k \in \pnk$, with $\mathcal{P}_i \in \pr^{n_i}$,
corresponds to a prime $\N^k$-homogeneous ideal $I_P$
of height $\sum_{i=1}^k n_i$ in $R$.  Furthermore,
$I_P = (L_{1,1},\ldots,L_{1,n_1},\ldots,L_{k,1},\ldots,L_{k,n_k})$
where $\deg L_{i,j} = e_i$ and $(L_{i,1},\ldots,L_{i,n_i})$
is the defining ideal of $\mathcal{P}_i \in \pr^{n_i}$.  If $\X 
= \{P_1,\ldots,P_s\} \subseteq \pnk$, then $\Ix = \bigcap_{i=1}^s
I_{P_i}$, where $I_{P_i}$ corresponds to $P_i$, is
the $\N^k$-homogeneous ideal of $R$ associated to $\X$.  The ring $R/\Ix$
inherits an $\N^k$-graded structure.  The Hilbert
function of $\X$ is then defined by $H_{\X}(\bi) =
\dim_{\bf k} (R/\Ix)_{\bi}$ for all $\bi \in \N^k$.  In this paper 
we study these Hilbert functions,
thereby building upon \cite{GuMaRa1,VT2}.

Each ideal $I_{P_i}$ is also homogeneous with respect to
the standard grading, so $I_{P_i}$ defines a linear variety
of dimension $k-1$ in $\pr^{N-1}$ where $N = \sum_{i=1}^k (n_i+1)$.
One can therefore take the point of view that our 
investigation of sets of points in $\pnk$ is an investigation of 
reduced unions of linear varieties  with extra 
hypotheses on the generators to ensure $\Ix$ is $\N^k$-homogeneous.

Ideally, we would like to 
classify those functions that arise as the Hilbert function
of a  set of points in $\pnk$.  However, besides the case  $k = 1$
which is dealt with in \cite{GeGrRo} and \cite{GeMaRo}, 
such a classification continues
to elude us. 
Though some properties of the Hilbert function
are known if $k >1$ (cf. \cite{GuMaRa1,VT2}),
even for sets of points in $\popo$
this problem remains open.

The proof of the characterization for the case $k = 1$
relies, in part, on the fact that the coordinate
ring of any finite set of points in $\pr^n$ is always Cohen-Macaulay (CM).
However, if $k > 1$, we show how to construct sets of points
which fail to be CM.  In fact, for each integer $l \in \{1,\ldots,k\}$,
we can construct a set of points with $\operatorname{depth} R/\Ix
= l$.  The failure of $R/\Ix$ to be CM in general provides an obstruction
to generalizing the proofs of \cite{GeGrRo,GeMaRo}.

We therefore restrict our investigation to sets of points
that arithmetically Cohen-Macaulay (ACM).
With this extra hypothesis on our set of points, we can generalize the proof 
for the case $k =1$
as given in \cite{GeGrRo} to all $k$.  In particular, we show that $H_{\X}$
is the Hilbert function of an ACM set of points in $\pnk$
if and only if $\Delta H_{\X}$, a generalized first difference function,
is the Hilbert function of some $\N^k$-graded artinian quotient.
Our generalization relies on two main ingredients: (1)
the existence of a regular sequence in $R/\Ix$ such that
each element has a specific multi-degree, and (2) the techniques
of \cite{MN} for lifting monomial ideals.

This characterization is not very satisfactory because it
translates our original problem into the open problem
of characterizing the Hilbert functions of $\N^k$-graded artinian
quotients.  However, we  characterize these quotients
in the special case $n_1 = \cdots = n_k = 1$, thereby giving
a complete description of the Hilbert functions of ACM sets
of points in $\pr^1 \times \cdots \times \pr^1$. 

In the last two sections  we specialize to ACM
sets of points in $\popo$.  It was first shown in \cite{GuMaRa1} 
that the ACM sets of points are characterized by their Hilbert
function.  We give a new proof of this result, plus a new characterization
that depends only upon numerical information describing the 
set $\X$.  We then show that
this numerical information also enables us to completely
calculate the graded Betti numbers of the minimal
free resolution of $\Ix$ (and thus, $H_{\X}$) 
provided $\X$ is ACM.  This generalizes the fact
that the Hilbert function and Betti numbers of a set of $s$ points in
$\pr^1$ depend only upon $s$.


\section{Preliminaries: Multi-graded rings, Hilbert functions, Points}

Throughout this paper, {\bf k} denotes an algebraically
closed field of characteristic zero.
In this section we provide the necessary facts and definitions
about multi-graded rings, Hilbert functions,  and sets
of points in multi-projective spaces.  See
\cite{Ro,VT2} for more on these topics.

Let $\N:=\{0,1,2,\ldots\}$.  
For an integer $n \in \N$, we set $[n]:=\{1,\ldots,n\}$.
We denote $(i_1,\ldots,i_k)\in \N^k$ by $\bi$.  
We set $|\bi| := \sum_h i_h$.
If $\bi,\bj \in \N^k$, then $\bi + \bj := (i_1 + j_1, \ldots,i_k + j_k)$.
We write $\bi \geq \bj$ if $i_h \geq j_h$ for every $h= 1,\ldots,k$.  
The set $\N^k$ is a semi-group generated by $\{e_1,\ldots,e_k\}$
where $e_i := (0,\ldots,1,\ldots,0)$ is the $i^{th}$ standard basis vector 
of $\N^k$.  If $c \in \N$, then 
$ce_i:= (0,\ldots,c,\ldots,0)$ with $c$ in the $i^{th}$ position.

Set  $R = {\bf k}[x_{1,0},\ldots,x_{1,n_1},x_{2,0},\ldots,x_{2,n_2},
\ldots,x_{k,0},\ldots,x_{k,n_k}]$, and induce an $\N^k$-grading on $R$
by setting $\deg x_{i,j} = e_i$.  An element $x \in R$ is said to
be $\N^k$-{\it homogeneous} (or simply {\it homogeneous}
if the grading is clear) if $x \in R_{\bi}$
for some $\bi \in \N^k$.  If $x$ is homogeneous, then
$\deg x := \bi$. 
An ideal $I = (F_1,\ldots,F_r) \subseteq R$ is an  $\N^k$-{\it homogeneous}
(or simply, {\it homogeneous}) {\it ideal} 
if each $F_j$ is $\N^k$-homogeneous.

For every $\bi \in \N^k$, the set  $R_{\bi}$ is a finite dimensional vector 
space over ${\bf k}$.  Since a basis for $R_{\bi}$ is the set of 
all monomials of degree  $\bi$,
$\dim_{\bf k} R_{\bi} = 
\binom{n_1+i_1}{i_1}\binom{n_2 + i_2}{i_2}\cdots \binom{n_k + i_k}{i_k}$.

If $I \subseteq R$ is a homogeneous ideal, then 
$S = R/I$ inherits an
$\N^k$-graded ring structure if we define $S_{\bi} = (R/I)_{\bi} :=
R_{\bi}/I_{\bi}$.
The numerical function $H_S:\N^k \rightarrow \N$ defined by
$H_S(\bi) := \dim_{\bf k} (R/I)_{\bi} 
= \dim_{\bf k} R_{\bi} - \dim_{\bf k} I_{\bi}$
is the {\it Hilbert function of $S$}.  
If $H:\N^k \rightarrow \N$ is a numerical function,
then the
{\it first difference function of $H$}, denoted $\Delta H$, is defined by 
\[
\Delta H(\bi) := \sum_{\underline{0} \leq \underline{l} = 
(l_1,\ldots,l_k) \leq (1,\ldots,1)} (-1)^{|\underline{l}|} H(i_1-l_1,\ldots,
i_k-l_k),
\]
where $H(\bj) = 0$ if $ \bj \not\geq \underline{0}$.  If
$k =1$, then our definition agrees with the classical definition.

The $\N^k$-graded polynomial ring $R$ is the coordinate ring of $\pnk$.
Let 
\[
P = [a_{1,0}:\ldots:a_{1,n_1}]\times \cdots \times 
[a_{k,0}:\ldots:a_{k,n_k}] \in \pnk
\]
be a point in this space.
The ideal of $R$ associated to the point $P$ is the prime ideal
$I_P = (L_{1,1},\ldots,L_{1,n_1},\ldots,L_{k,1},\ldots,L_{k,n_k})$
where $\deg L_{i,j} = e_i$ for $j = 1,\ldots,n_i$.
If $P_1,\ldots,P_s$ are $s$ distinct points
and $\X = \{P_1,\ldots,P_s\} \subseteq \pnk$,
then $\Ix = I_{P_1} \cap \cdots \cap I_{P_s}$ where
$I_{P_i}$ is the ideal associated to the point $P_i$.
The ring $R/\Ix$ then has the following property:

\begin{lemma}[{\cite[Lemma 3.3]{VT}}] \label{firstnzd}
Let $\X$ be a finite set of  points in $\pnk$.
For each integer $i \in [k]$  there exists
a form $L_i \in R_{e_i}$ such that $\overline{L}_i$
is a non-zero divisor in $R/\Ix$.  
\end{lemma}

We write $H_{\X}$ to denote the Hilbert function $H_{R/\Ix}$,
and we  say $H_{\X}$ is the {\it Hilbert function}
of $\X$.  Let $\pi_i:\pnk \rightarrow \pr^{n_i}$ denote the $i^{th}$
projection morphism.  Then $t_i:=|\pi_i(\X)|$ is the
number of distinct $i^{th}$-coordinates in $\X$. 
With this notation we have

\begin{proposition} \label{eventualgrowth}
Let $\X$ be a finite set of points in $\pnk$ with Hilbert function 
$H_{\X}$.   
\begin{enumerate}
\item[$(i)$] {\cite[Proposition 4.2]{VT2}}
If $(j_1,\ldots,j_i,\ldots,j_k) \in \N^k$ and if $j_i \geq t_i - 1,$ then 
\[ 
H_{\X}(j_1,\ldots,j_i,\ldots,j_k) = H_{\X}(j_1,\ldots,t_i-1,\ldots,j_k).
\]
\item[$(ii)$] {\cite[Corollary 4.8]{VT2}} 
$H_{\X}(j_1,\ldots,j_k) = s $ for all
$(j_1,\ldots,j_k) \geq (t_1-1,t_2-1,\ldots,t_k-1).$  
\end{enumerate}
\end{proposition}

\begin{remark} 
Fix an integer $i \in [k]$, and fix $k-1$ integers in $\N$,
say $j_1,\ldots,j_{i-1},j_{i+1},\ldots,j_k$.  Set
$\underline{j}_{l} := (j_1,\ldots,j_{i-1},l,j_{i+1},\ldots,j_k)$
for each integer $l \in \N$.  Then
Proposition \ref{eventualgrowth} (i) can be interpreted as saying the sequence
$\left\{H_{\X}(\underline{j}_l)\right\}$
becomes constant.  In fact, $H_{\X}(\underline{j}_l)
= H_{\X}(\underline{j}_{t_i-1})$ for all $l \geq t_i -1$.
\end{remark}

\begin{proposition}[{\cite[Proposition 3.2]{VT2}}] \label{sides}
Let $\X$ be a finite set of points in $\pnk$ with
Hilbert function $H_{\X}$.  Fix an integer $i \in [k]$.
Then the sequence $H = \{h_j\}$, where $h_j := H_{\X}(je_i)$, 
is the Hilbert function of $\pi_i(\X) \subseteq \pr^{n_i}$.
\end{proposition}

We end this section with some comments on the depth
and Krull dimension of $R/\Ix$.
Let ${\bf m} := \bigoplus_{\underline{0} \neq \underline{j} \in \N^k} R_{\bj}
= (x_{1,0},\ldots,x_{k,n_k})$
be the maximal ideal of $R$.  
If $I \subseteq R$ is an $\N^k$-homogeneous ideal, then 
recall that we say a sequence $F_1,\ldots,F_r$ of 
elements is a {\it regular sequence modulo $I$} if and only if
$(i)$  $(I,F_1,\ldots,F_r) \subseteq {\bf m}$,
$(ii)$ $\overline{F}_1$ is not a zero divisor in $R/I$, and 
$(iii)$ $\overline{F}_i$ is not a zero divisor in
$R/(I,F_1,\ldots,F_{i-1})$ for $1 < i \leq r$.
The {\it depth of $R/I$}, written $\dep R/I$, is the length
of the maximal regular sequence modulo $I$.

Because each prime ideal $I_{P_i}$ has height $\sum_{i=1}^k n_i$,
it follows that K-$\dim R/\Ix = k$, the number of projective
spaces.  This result, coupled with Lemma \ref{firstnzd},
implies $1 \leq \operatorname{depth} R/\Ix \leq k$.
Thus,  every set of points in $\pr^n$ has $\dep R/\Ix = 1$.  
If $k \geq 2$, the value for $\dep R/\Ix$ is not immediately clear.  
In fact, for each integer $l \in [k]$
we can construct a set
of points in $\X$ such that $\dep R/\Ix = l$.
We begin with a lemma.

\begin{lemma} \label{depthlemma}
Fix a positive integer $k$.  Denote by  $X_1$ and $X_2$ the two points
\[
X_1 := [1:0] \times [1:0] \times \cdots \times [1:0] ~\text{and}~~
X_2 := [0:1] \times [0:1] \times \cdots \times [0:1]
\]
in $\pr^1\times\cdots\times\pr^1$ (k times).
If $\X := \{X_1,X_2\}$,  then $\operatorname{depth} R/\Ix = 1.$
\end{lemma}

\begin{proof}
The defining ideal of $\X$ is
$\Ix = I_{X_1} \cap I_{X_2} =  \left( \{x_{a,0}x_{b,1} \left|~ 1\leq a \leq k, 
1\leq b \leq k\}\right)\right.$
in the $\N^k$-graded ring 
$R = {\bf k}[x_{1,0},x_{1,1},x_{2,0},x_{2,1},\ldots,x_{k,0},x_{k,1}]$.
Since $x_{1,0} + x_{1,1}$ does not vanish at either 
point, it suffices to show that
every non-zero element of $R/(\Ix,x_{1,0}+x_{1,1})$ is a zero divisor.

So, set $J =(\Ix,x_{1,0}+x_{1,1})$ and  suppose that $0 \neq
\overline{F} \in R/J$.  Without loss of generality,
we can take $F$ to be $\N^k$-homogeneous.  We write $F$ as 
$F = F_0 + F_1x_{1,0} + F_2x_{1,0}^2 + \cdots$
where $F_i \in {\bf k}[x_{1,1},x_{2,0},\ldots,x_{k,1}]$.
Since $x_{1,0}x_{1,1} \in \Ix$, it follows that 
$x_{1,0}^2 = x_{1,0}(x_{1,0} + x_{1,1}) - x_{1,0}x_{1,1} \in J$.
Hence, we can assume that $F = F_0 + F_1x_{1,0}.$
The element $x_{1,0} \not\in J$.  Although $x_{1,0}, F\not\in J$,
we claim that $Fx_{1,0} \in J$.  Indeed, for each integer
$1 \leq b \leq k$,  $x_{1,0}x_{b,1} \in I_{\X} \subseteq J$.  Furthermore,
for  $1 \leq a \leq k$, the element 
$x_{1,0}x_{a,0} = x_{a,0}\left(x_{1,0} + x_{1,1}\right) - x_{a,0}x_{1,1}
\in J$.  Hence, each term of $F_0x_{1,0}$
is in $J$, so $F_0x_{1,0} \in J$.
But then $Fx_{1,0} = F_0x_{1,0} + F_1x_{1,0}^2 \in
J$.  So, every $0 \neq \overline{F} \in R/J$ is a zero divisor because it
is annihilated by  $\overline{x}_{1,0}$.
\end{proof}

\begin{proposition} \label{possibledepths}
Fix a positive integer $k$, and let $n_1,\ldots,n_k$
be  $k$ positive integers.  Then, for every integer $l
\in [k]$, there exists
a set  of points $\X$ in $\pnk$ with 
$\operatorname{depth} R/\Ix = l$.
\end{proposition}

\begin{proof}
For every $l \in [k]$ we construct a set  with
the desired depth.  Set $P_i := [1:0:\cdots:0] \in \pr^{n_i}$
for $1 \leq i \leq k$ and $Q_i := [0:1:0\cdots:0] \in \pr^{n_i}$
for $1 \leq i \leq k$.
Fix an $l \in [k]$ and let $X_1$ and $X_2$ be
the following two points of $\pnk$:
\[
X_1 := P_1 \times P_2 \times \cdots \times P_k ~~~\text{and}~~~
X_2 := P_1 \times P_2 \times \cdots \times P_{l-1} 
\times Q_l \times \ldots \times Q_k.
\]
If we set  $\X_l := \{X_1,X_2\}$, we claim that 
$\operatorname{depth} R/I_{\X_l}=l$.
The defining ideal of $\X_l$ is
\begin{eqnarray*}
I_{\X_l} & =  &
 \left(
\begin{array}{l}
x_{1,1},\ldots,x_{1,n_1},\ldots,x_{l-1,1},\ldots,x_{l-1,n_{l-1}},\\
x_{l,2},\ldots,x_{l,n_l},\ldots,x_{k,2},\ldots,x_{k,n_k}, \\
\{x_{a,0}x_{b,1} ~|~ l \leq a \leq k,~ l \leq b \leq k\}
\end{array}\right).
\end{eqnarray*}
It follows that $R/I_{\X_l} \cong S/J$ where 
\[ 
S/J =  \frac{{\bf k}[x_{1,0},x_{2,0},x_{3,0},\ldots,x_{l-1,0},
x_{l,0},x_{l,1},x_{l+1,0},x_{l+1,1},\ldots,x_{k,0},x_{k,1}]}
{(\{x_{a,0}x_{b,1} ~|~ l \leq a \leq k, l \leq b \leq k\})}.
\]
The indeterminates $x_{1,0},x_{2,0},\ldots,x_{l-1,0}$
give rise to a regular sequence in $S/J$.  Thus, 
$\operatorname{depth} R/I_{\X_l} = \operatorname{depth} S/J \geq l-1.$
Set $K = (J,x_{1,0},\ldots,x_{l-1,0})$.  Then
\[
S/K \cong 
\frac{{\bf k}[x_{l,0},x_{l,1},x_{l+1,0},x_{l+1,1},\ldots,x_{k,0},x_{k,1}]}
{(\{x_{a,0}x_{b,1} ~|~ l \leq a \leq k, l \leq b \leq k\})}.
\]
The ring  $S/K$ is then isomorphic to the  $\N^{k-l+1}$-graded
coordinate ring of the  two points 
$\left\{[1:0]\times [1:0] \times \cdots \times [1:0],
[0:1] \times [0:1]\times\cdots \times [0:1]\right\}$
in ${\pr^1\times \cdots \times \pr^1}$ ($(k-l+1)$ times).  From  
Lemma \ref{depthlemma} we have  $\operatorname{depth} S/K = 1$, and hence,
$\operatorname{depth} R/I_{\X_l} = l-1 + 1 = l.$
\end{proof}


\section{The Hilbert functions of ACM Sets of Points}

For an arbitrary set of points $\X \subseteq \pnk$,  
$\operatorname{depth} R/\Ix \leq  $ K-$\dim R/\Ix$.  
If the equality holds, the coordinate ring 
is Cohen-Macaulay (CM), and the set of points
are said to be arithmetically Cohen-Macaulay (ACM).  
We now investigate the Hilbert functions
of those sets of points in $\pnk$ that are also ACM.  Under
this extra hypothesis we can generalize
the characterization of the Hilbert functions of
sets of points in $\pr^n$ as found in \cite{GeGrRo}.

We begin with a preparatory lemma.

\begin{lemma}     	\label{idealelements}
Let $\X \subseteq \pnk$ be a finite set of points, and 
suppose $L_1,\ldots,L_t$, with $t \leq k$ and $\deg L_i = e_i$,
give rise to a regular sequence in $R/\Ix$.  Then there
exists a positive integer $l$ such that
$(x_{1,0},\ldots,x_{1,n_1},\ldots,x_{t,0},\ldots,x_{t,n_t})^l 
\subseteq (\Ix,L_1,\ldots,L_{t}).$
\end{lemma}

\begin{proof}
Set $J_i := (\Ix,L_1,\ldots,L_i)$ for $i = 1,\ldots,t$.
Since $L_1,\ldots,L_t$ form a regular sequence on $R/\Ix$,
for each $i = 1,\ldots,t$ we have a short exact sequence with 
degree $(0,\ldots,0)$ maps:
\[
0 \rightarrow \left(R/J_{i-1}\right)(-e_i) 
\stackrel{\times \overline{L}_i}{\rightarrow}
R/J_{i-1} \rightarrow R/J_i \rightarrow 0
\]
where we set  $J_0 :=\Ix$.
From the  exact sequences we derive the following formula:
\[
\dim_{\bf k} (R/J_t)_{\bi}
= \sum_{\underline{0} \leq (j_1,\ldots,j_t) \leq \underline{1}}
(-1)^{(j_1+\cdots+j_t)} \dim_{\bf k} (R/\Ix)_{i_1-j_1,\ldots,i_t-j_t,
i_{t+1},\ldots,i_k}
\]
where we set $\dim_{\bf k} (R/\Ix)_{\bj} = 0$ if $\bj \not\geq \underline{0}$.

For each integer $j = 1,\ldots,t$, set $t_j := |\pi_j(\X)|$.
By Proposition \ref{eventualgrowth}, if $i_j \geq t_j$, then
$\dim_{\bf k}(R/\Ix)_{i_je_j} =  \dim_{\bf k} (R/\Ix)_{(i_j-1)e_j}.$
This fact, coupled with above above formula, implies that $
\dim_{\bf k} (R/J_t)_{t_je_j}  = 
\dim_{\bf k} (R/\Ix)_{t_je_j} -
\dim_{\bf k} (R/\Ix)_{(t_j-1)e_j} =  0.
$
Thus $R_{t_je_j}
= (J_t)_{t_je_j}$, or equivalently,
$(x_{j,0},\ldots,x_{j,n_j})^{t_j} \subseteq
(\Ix,L_1,\ldots,L_t)$.  Since this is true for
each integer $1 \leq j \leq t$, there exists an integer $l \gg 0$
such that 
$
(x_{1,0},\ldots,x_{1,n_1},\ldots,x_{t,0},\ldots,x_{t,n_t})^l 
\subseteq J_t$,
as desired.
\end{proof}

\begin{proposition}	\label{regularsequence}
Suppose that $\X$ is an ACM set of points in $\pnk$.
Then there exist elements $\overline{L}_1,\ldots,\overline{L}_k$
in $R/\Ix$ such that $L_1,\ldots,L_k$ give rise to a regular
sequence in $R/\Ix$ and $\deg L_i = e_i$.
\end{proposition}

\begin{proof}
The existence of a regular sequence of length $k$ follows
from the definition of a CM ring.  The non-trivial
part of this statement is the existence of a regular
sequence whose elements have specific multi-degrees.

By Lemma \ref{firstnzd} there exists a form $L_1 \in R_{e_1}$
such that $\overline{L}_1$ is a non-zero divisor of $R/\Ix$.
To complete the proof it is enough to show 
for each $t = 2,...,k$ there exists
an element $L_t \in R_{e_t}$ such that $\overline{L}_t$
is a non-zero divisor of the ring $R/(\Ix,L_1,\ldots,L_{t-1})$.

Set $J := (\Ix,L_1,\ldots,L_{t-1})$ and  let $J = Q_1 \cap \cdots \cap Q_r$
be the ideal's primary decomposition.  
For each $i$ set $\wp_i := \sqrt{Q_i}$.  
Since  $\bigcup_{i=1}^r \overline{\wp}_i$ is
the set of zero divisors of $R/J$,
we want to show that $\bigcup_{i=1}^r (\wp_i)_{e_t}
\subsetneq R_{e_t}$.  Since $R_{e_t}$ is a vector space over
an infinite field, $R_{e_t}$ cannot be expressed as
the finite union of proper subvector spaces.  Thus, it
suffices to show $(\wp_i)_{e_t} \subsetneq R_{e_t}$
for each $i$.

So, suppose there is $i \in [r]$ such that
$(\wp_i)_{e_t} = R_{e_t}$, or equivalently, 
$(x_{t,0},\ldots,x_{t,n_t}) \subseteq \wp_i$.  
By Lemma \ref{idealelements} there exists 
$l \in \N^+$ such that
$(x_{1,0},\ldots,x_{1,n_1},\ldots,x_{t-1,0},\ldots,x_{t-1,n_{t-1}})^l
\subseteq J \subseteq Q_i$.  It follows that
$(x_{1,0},\ldots,x_{1,n_1},\ldots,x_{t,0},\ldots,x_{t,n_t}) \subseteq \wp_i.$

Because the prime ideal $\wp_i$ also contains $\Ix = I_{P_1}
\cap \cdots \cap I_{P_s}$, where $I_{P_j}$ is the prime
ideal associated to $P_j \in \X$, we
can assume, after relabeling, $I_{P_1} \subseteq \wp_i$.  Let $\wp := I_{P_1}
+(x_{1,0},\ldots,x_{t,n_t})$.  
Since
$I_{P_1} = (L_{1,1},\ldots,L_{1,n_1},\ldots,L_{k,1},\ldots,L_{k,n_k})$
where $\deg L_{m,n} = e_m$
\[
\wp = (x_{1,0},\ldots,x_{t,n_t},L_{t+1,1},\ldots,L_{t+1,n_{t+1}},
\ldots,L_{k,1},\ldots,L_{k,n_k}) \subseteq \wp_i.
\]
Thus $\operatorname{ht}_R(\wp_i) \geq  
\operatorname{ht}_R(\wp) = \left(\sum_{i=1}^k n_i\right) +t$,
where $\operatorname{ht}_R(I)$ denotes the height of $I$.

From the identity $\operatorname{ht}_R (J)  =  \kdim R - \kdim R/J$
we calculate the height of $J$:
\[
\operatorname{ht}_R (J)  =  \left(\sum_{i=1}^k n_i +1 \right) - (k - (t-1))
=  \left(\sum_{i=1}^k n_i \right) + (t-1).
\]
Since $\X$ is ACM, $R/J$ is CM, and hence the ideal $J$ is height
unmixed, i.e., all the associated primes of $J$
have height equal to $\operatorname{ht}_R(J)$.  But $\wp_i$ is an 
associated prime of $J$ with 
$
\operatorname{ht}_R (\wp_i) > \operatorname{ht}_R (J).
$
This contradiction implies our assumption
$(\wp_i)_{e_t} = R_{e_t}$ cannot be true.
\end{proof}

\begin{remark} If $S = {\bf k}[x_0,\ldots,x_n]$
is an $\N^1$-graded ring with $I \subseteq S$ such that $S/I$ is CM, 
then a maximal regular sequence can be chosen so that each
element is homogeneous \cite[Proposition 1.5.11]{BrHe}.
However, as stated in \cite{St} (but no example is given),
it is not always possible to pick a regular sequence that
respects the multi-grading.  For example, let $S = {\bf k}[x,y]$ with
$\deg x = (1,0)$ and $\deg y = (0,1)$ and $I = (xy)$.  Then $S/I$
is CM, but all homogeneous elements of $S/I$,
which  have the form
$\overline{cx}^a$ or $\overline{cy}^b$ with $c \in {\bf k}$,
are  zero divisors.  Note that $\overline{x} 
+ \overline{y}$ is a non-zero divisor, but not homogeneous.
The fact that a homogeneous regular sequence can be found
in a multi-graded ring is thus a very special situation.
\end{remark}

We extend the notion of a graded artinian quotient
in the natural way.

\begin{definition}
A homogeneous ideal $I$ in the $\N^k$-graded ring  
$R$  is an {\it artinian ideal} if any of the following equivalent
statements hold:
\begin{enumerate}
\item[$(i)$]  $\operatorname{K-}\dim R/I =0$.
\item[$(ii)$] $\sqrt{I} = {\bf m} =(x_{1,0},\ldots,x_{1,n_1},
\ldots,x_{k,0},\ldots,x_{k,n_k})$.
\item[$(iii)$]  For each integer $i \in [k]$, 
$H_{R/I}(le_i) = 0$ for all $l \gg 0$.
\end{enumerate}
A ring $S = R/I$ is an {\it $\N^k$-graded artinian quotient} 
if the homogeneous ideal $I$ is an artinian. 
\end{definition}

\begin{corollary}	\label{ahilbert2}
Let $\X$ be an ACM set of points in $\pnk$
with Hilbert function $H_{\X}$.  Then
\[
\Delta H_{\X}(i_1,\ldots,i_k) := \sum_{\underline{0} \leq
\underline{l} = 
(l_1,\ldots,l_k) \leq (1,\ldots,1)} (-1)^{|\underline{l}|} 
H_{\X}(i_1-l_1,\ldots,i_k-l_k),
\]
where $H_{\X}(\bi) = 0$ if $\bi \not\geq 
\underline{0}$, is the Hilbert function of some $\N^k$-graded artinian 
quotient of the ring ${\bf k}[x_{1,1},\ldots,x_{1,n_1},\ldots,
x_{k,1},\ldots,x_{k,n_k}]$.
\end{corollary}

\begin{proof}  
By Proposition \ref{regularsequence} there exists $k$ forms
$L_1,\ldots,L_k$ that give rise to a regular sequence
in $R/\Ix$ with $\deg L_i = e_i$.  After
making a linear change of variables in the
$x_{1,i}$'s, a linear change of variables in the
$x_{2,i}$'s, etc., we can assume that
$L_i = x_{i,0}$.

The ideal  $(\Ix,x_{1,0},\ldots,x_{k,0})/(x_{1,0},\ldots,x_{k,0})$
is isomorphic to an ideal $J$ of the  ring
$S = {\bf k}[x_{1,1},\ldots,x_{1,n_1},\ldots,
x_{k,1},\ldots,x_{k,n_k}]$.  Set $A := S/J$, and so
\[
A \cong \frac{R/(x_{1,0},\ldots,x_{k,0})}
{(\Ix,x_{1,0},\ldots,x_{k,0})/(x_{1,0},\ldots,x_{k,0})}
\cong  \frac {R}{(\Ix,x_{1,0},\ldots,x_{k,0})}.
\]
The ring $A$ is artinian because  
there exists $l \gg 0$ by Lemma \ref{idealelements} such that 
${\bf m}^l \subseteq
(\Ix,x_{1,0},\ldots,x_{k,0})$.

It therefore remains to compute
the Hilbert function of $A$. 
Set $J_i = (\Ix,x_{1,0},\ldots,x_{i,0})$ for $i = 1,\ldots,k$.
For each $i = 1,\ldots,k$ we have a short exact sequence with 
degree $(0,\ldots,0)$ maps:
\[
0 \longrightarrow \left(R/J_{i-1}\right)(-e_i)
\stackrel{\times \overline{x}_{i,0}}{\longrightarrow}
R/J_{i-1} \longrightarrow R/J_i \longrightarrow 0
\]
where $J_0 := \Ix$.
From the $k$ short exact sequences we have that
\[H_{R/J_k}(\bi) = \Delta H_{\X}(\bi) 
 := \sum_{\underline{0} \leq \underline{l} = 
(l_1,\ldots,l_k) \leq (1,\ldots,1)} (-1)^{|\underline{l}|} 
H_{\X}(i_1-l_1,\ldots,i_k-l_k)\]
where $H_{\X}(\bi) = 0$ if $\bi \not\geq \underline{0}$.
This completes the proof since $A \cong R/J_k$.
\end{proof}

The remainder of this section is devoted to showing
the necessary condition in 
Corollary ~\ref{ahilbert2} is also sufficient.
To demonstrate this converse, 
we describe how to {\it lift} an ideal.

\begin{definition}  	\label{lifting}
Let $R = {\bf k}[x_{1,0},\ldots,x_{1,n_1},\ldots,x_{k,0},\ldots,x_{k,n_k}]$ 
and let $S= {\bf k}[x_{1,1},\ldots,
x_{1,n_1},\ldots,\newline
x_{k,1},\ldots,x_{k,n_k}]$ be $\N^k$-graded rings.  Let 
$I \subseteq R$ and $J \subseteq S$ be $\N^k$-homogeneous ideals.  Then we say
$I$ is a {\it lifting of $J$ to $R$} if
\begin{enumerate}
\item[$(i)$]  $I$ is radical in $R$
\item[$(ii)$] ${\displaystyle 
\frac{(I,x_{1,0},\ldots,x_{k,0})}{(x_{1,0},\ldots,x_{k,0})} \cong J}$
\item[$(iii)$] $x_{1,0},\ldots,x_{k,0}$ give rise to a regular sequence in
$R/I$.
\end{enumerate}
\end{definition}

Our plan is to
lift a monomial ideal of $S$ to an $\N^k$-homogeneous
ideal $I$  of $R$, using the techniques and results
of \cite{MN}, so such that $I$ is the ideal of a reduced
set of points in $\pnk$. 
We make a brief digression to introduce the relevant content of \cite{MN}.

Suppose that $S$ and $R$ are
as in Definition \ref{lifting}, but for the moment, we
only assume that they are $\N^1$-graded.
 For each indeterminate
$x_{i,j}$ with $1\leq i \leq k$ and $1\leq j \leq n_i$,
choose infinitely many linear forms $L_{i,j,l} \in 
{\bf k}[x_{i,j},x_{1,0},x_{2,0},\ldots,x_{k,0}]$ with $l\in \N^+$.
We only assume that the coefficient of $x_{i,j}$ in $L_{i,j,l}$ is not
zero.  The infinite matrix $A$, where 
\[ A:= \bmatrix
L_{1,1,1} & L_{1,1,2} & L_{1,1,3} & \cdots \\
\vdots & \vdots & \vdots & \\
L_{1,n_1,1} & L_{1,n_1,2} & L_{1,n_1,3} & \cdots \\
\vdots & \vdots & \vdots & \\
L_{k,1,1} & L_{k,1,2} & L_{k,1,3} & \cdots \\
\vdots & \vdots & \vdots & \\
L_{k,n_k,1} & L_{k,n_k1,2} & L_{k,n_k,3} & \cdots \\
\endbmatrix\]
is called a {\it lifting matrix}.
By using the lifting matrix, we associate to
each monomial $m = x_{1,1}^{a_{1,1}}\cdots x_{1,n_1}^{a_{1,n_1}}\cdots
x_{k,1}^{a_{k,1}}\cdots x_{k,n_k}^{a_{k,n_k}}$ of the ring $S$ the
element
\[\overline{m} = 
\left[ \prod_{i=1}^{n_1}\left(\prod_{j=1}^{a_{1,i}} L_{1,i,j}
\right)\right]\cdots
\left[ \prod_{i=1}^{n_k}\left(\prod_{j=1}^{a_{k,i}}
L_{k,i,j}\right)\right]
\in R.\]
Depending upon our choice of $L_{i,j,l}$'s, $\overline{m}$
may or may not be $\N^k$-homogeneous.  However, $\overline{m}$
is homogeneous.  If $J = (m_1,\ldots,m_r)$ is a monomial
ideal of $S$, then we use $I$ to denote the ideal
$(\overline{m}_1,\ldots,\overline{m}_r) \subseteq R$. 
The following properties, among others, relate $R/I$ and $S/J$.

\begin{proposition}[{\cite[Corollary 2.10]{MN}}] \label{liftingproperties}
Let $J \subseteq S$ be a monomial ideal, and let $I$ be the ideal
constructed from $J$ using any lifting matrix.
Then
\begin{enumerate}
\item[$(i)$] $S/J$ is CM if and only if $R/I$ is CM.
\item[$(ii)$]  ${\displaystyle 
\frac{(I,x_{1,0},\ldots,x_{k,0})}{(x_{1,0},\ldots,x_{k,0})} 
\cong J}$
\item[$(iii)$]  $x_{1,0},\ldots,x_{k,0}$ give rise to a regular sequence in
$R/I$.
\end{enumerate}
\end{proposition}

We now consider the lifting of a monomial ideal using the lifting matrix
$\mathcal{A} := [L_{i,j,l}]$ where
\[L_{i,j,l} = x_{i,j} - (l-1)x_{i,0} ~~~\mbox{for 
$1 \leq i \leq n, ~1\leq j \leq n_i$, and $l \in \N^+$.}\]
The lifting matrix $\mathcal{A}$ associates to every monomial 
$m =  x_{1,1}^{a_{1,1}}\cdots x_{1,n_1}^{a_{1,n_1}}\cdots
x_{k,1}^{a_{k,1}}\cdots x_{k,n_k}^{a_{k,n_k}}$ of $S$ the 
following $\N^k$-homogeneous form of $R$:
\[
\overline{m} = 
\left[ \prod_{i=1}^{n_1}\left(\prod_{j=1}^{a_{1,i}} (x_{1,i} - (j-1)x_{1,0})
\right)\right]\cdots
\left[ \prod_{i=1}^{n_k}\left(\prod_{j=1}^{a_{k,i}}
(x_{k,i} - (j-1)x_{k,0})\right)\right]
.\]
Thus, if $I$ is constructed from a monomial ideal $J \subseteq S$ using 
$\mathcal{A}$, $I$ is $\N^k$-homogeneous.
In fact,  using $\mathcal{A}$,
the ideal $I$ is a lifting of $J$ to $R$.  To prove this
statement, we require

\begin{lemma}[{\cite[Corollary 2.18]{MN}}] \label{primarylift}
Let $J \subseteq S$ be a monomial ideal and let $I$ be
be constructed from $J$ using any lifting matrix. If
$J = Q_1 \cap \cdots \cap Q_r$ is the primary
decomposition of $J$, then
$I = \overline{Q}_1 \cap \cdots \cap \overline{Q}_r$
where $\overline{Q}_i$ is the ideal generated by the lifting
of the generators of $Q_i$.
\end{lemma}

\begin{proposition} \label{liftingofJ}
Let $J \subseteq S$ be a monomial ideal and let $I$ be the
ideal constructed from $J$ using the lifting matrix $\mathcal{A}$. 
Then $I$ is a lifting of $J$ to $R$.
\end{proposition}

\begin{proof}
Since Proposition \ref{liftingproperties} is true for any lifting matrix,
it suffices to show that $I$ is radical.  
Let $J = Q_1 \cap \cdots \cap Q_r$ be the primary decomposition
of $J$.  Since $J$ is a monomial ideal, then by Remark 2.19 of \cite{MN}
we have that each $Q_i$ is a complete intersection of the form
\[Q_i = \left(x_{1,i_{1,1}}^{a_{1,i_{1,1}}},\ldots,
x_{1,{i_{1,p_1}}}^{a_{1,i_{1,p_1}}},\ldots
x_{k,i_{k,1}}^{a_{k,i_{k,1}}},\ldots,
x_{k,{i_{k,p_k}}}^{a_{k,i_{k,p_k}}}\right)\]
with $a_{j,i_{j,l}} \geq 1$ for each variable that appears in $Q_i$.
Using the lifting matrix $\mathcal{A}$ we then have 
\[\overline{Q}_i =
\left(\prod_{l=1}^{a_{1,i_{1,1}}} L_{1,i_{1,1},l},
\ldots,\prod_{l=1}^{a_{1,i_{1,p_1}}} L_{1,i_{1,p_1},l},\ldots,
\prod_{l=1}^{a_{k,i_{k,1}}} L_{k,i_{k,1},l},
\ldots,\prod_{l=1}^{a_{k,i_{k,p_k}}} L_{k,i_{k,p_k},l}
 \right)\]
where $L_{i,j,l} = x_{i,j} - (l-1)x_{i,0}$.  But then $\overline{Q}_i$
is a reduced complete intersection.  It then follows from
Lemma \ref{primarylift} that $I$ must be radical.
\end{proof}

We now describe the zero set of the lifted ideal $I$.
For each
\[
\alphak := ((a_{1,1},\ldots,a_{1,n_1}),\ldots,(a_{k,1},\ldots,a_{k,n_k}))
\in \N^{n_1} \times \cdots \times \N^{n_k},
\]
set
$
\xka := x_{1,1}^{a_{1,1}}\cdots x_{1,n_1}^{a_{1,n_1}}\cdots
x_{k,1}^{a_{k,1}}\cdots x_{k,n_k}^{a_{k,n_k}}.
$  
If $P$ is the set of all monomials of $S$ including the monomial $1$,
then  there exists a bijection between $P$ and $\N^{n_1} \times \cdots 
\times \N^{n_k}$ given by the map $\xka \leftrightarrow \alphak$. 
To each tuple $\alphak$
we associate the point $\alphako \in \pnk$ where
\[
\alphako := [1:a_{1,1}:a_{1,2}:\cdots:a_{1,n_1}] 
\times \cdots \times [1: a_{k,1}: a_{k,2} : \cdots:a_{k,n_k}].
\]
Note that if $m = \xka \in P$ and if $\overline{m}$ is constructed from 
$\xka$ using the lifting matrix $\mathcal{A}$, 
then $\overline{m}(\alphako) \neq 0$. 
In fact, it follows  from our construction that 
$\overline{m}(\overline{(\underline{\beta}_1,\ldots,
\underline{\beta}_k)}) = 0$ if and only if
some coordinate of $(\underline{\beta}_1,\ldots,\underline{\beta}_k)$ 
is strictly less than some coordinate of $\alphak$.

If $J$ is a monomial ideal of $S$, then let $N$ be the set of monomials
in $J$.   The elements of $M := P \backslash N$
are representatives for a ${\bf k}$-basis of the $\N^k$-graded ring
$S/J$.  Set 
\[
\overline{M} := 
\left\{{\overline{(\underline{\beta}_1,\ldots,
\underline{\beta}_k)} \in \pnk} \left| {~X_1^{\underline{\beta}_1}
\cdots X_k^{\underline{\beta}_k} \in M} \right\}\right..
\]  
Let ${\bf I}(\overline{M})$ denote the $\N^k$-homogeneous ideal
associated to $\overline{M}$.  If 
$m = \xka \in J$ is a minimal generator, then for
each $X_1^{\underline{\beta}_1}
\cdots X_k^{\underline{\beta}_k} \in M$ there exists
at least one coordinate of $(\underline{\beta}_1,\ldots,
\underline{\beta}_k)$ that is strictly less then
some coordinate of $\alphak$.
So $\overline{m}(\overline{(\underline{\beta}_1,\ldots,
\underline{\beta}_k)}) = 0$ for all $X_1^{\underline{\beta}_1}
\cdots X_k^{\underline{\beta}_k} \in M$, and hence the
lifted ideal $I \subseteq {\bf I}(\overline{M})$.
On the other hand $\overline{M} = {\bf V}(I)$, the 
zero set of $I$, so by the $\N^k$-graded analog of the 
Nullstellensatz (cf. \cite[Theorem 2.3]{VT2}) we
have ${\bf I}(\overline{M}) \subseteq \sqrt{I}$.  Because
$I$ is radical by Proposition \ref{liftingofJ} we have just
shown

\begin{lemma}	\label{RADICAL1}
Let $J \subseteq S$ be a monomial ideal and let $I$ be the ideal
constructed from $J$ using the lifting matrix $\mathcal{A}$.
Then, with the notation as above, $I = {\bf I}(\overline{M})$.
\end{lemma}

We come to the main result of this section.

\begin{theorem}		\label{acmforpnk}
Let $H: \N^k \rightarrow \N$ be a numerical function.  
Then $H$ is the Hilbert function of an ACM  set of distinct points
in $\pnk$ if and only if the first difference function
\[
\Delta H(i_1,\ldots,i_k) =  \sum_{\underline{0} \leq \underline{l} = 
(l_1,\ldots,l_k) \leq (1,\ldots,1)} (-1)^{|\underline{l}|} 
H(i_1-l_1,\ldots,i_k-l_k),\] 
where $H(\bi) = 0$ if $\bi \not\geq \underline{0}$,
is the Hilbert function of some $\N^k$-graded artinian quotient of 
$S={\bf k}[x_{1,1},\ldots,x_{1,n_1},\ldots,x_{k,1},\ldots,x_{k,n_k}]$.
\end{theorem}

\begin{proof}
Because of Corollary \ref{ahilbert2}, we only need to show one direction.  
So, if $\Delta H$ is the Hilbert function
of some $\N^k$-graded artinian quotient of $S$, then there
exists an $\N^k$-homogeneous ideal $J \subseteq S$ with
$\Delta H(\bi) = H_{S/J}(\bi)$ for all $\bi \in \N^k$.  By replacing $J$ 
with its leading term ideal,  we can assume that $J = (m_1,\ldots,
m_r)$ is a monomial ideal of $S$.

Let $I \subseteq R$
be the ideal constructed from $J$ using the lifting matrix $\mathcal{A}$.  By
Proposition \ref{liftingofJ}, 
the ideal $J \cong (I,x_{1,0},\ldots,x_{k,0})/(x_{1,0},\ldots,x_{k,0}) $ 
where $x_{1,0},\ldots,x_{k,0}$ give rise to a regular sequence 
in $R/I$.  Because $\deg x_{i,0} = e_i$, we have $k$ short 
exact sequences with degree $(0,\ldots,0)$ maps:
\[
0 \longrightarrow \left(R/J_{i-1}\right)(e_i) 
\stackrel{\times \overline{x}_{i,0}}{\longrightarrow}
R/J_{i-1} \longrightarrow R/J_i \longrightarrow 0
\]
where $J_i := (I,x_{1,0},\ldots,x_{i,0})$ for $i = 1,\ldots,k$ and $J_0: = I$. 
Furthermore, 
\[
S/J \cong \frac{R/(x_{1,0},\ldots,x_{k,0})}{(I,x_{1,0},\ldots,x_{k,0})/(x_{1,0},\ldots,x_{k,0})}
\cong R/(I,x_{1,0},\ldots,x_{k,0}) = R/J_k.
\]
Then, using the $k$ short exact sequences to calculate the $H_{R/I}$,
we find that $H = H_{R/I}$.

 If $N$ is the set of monomials in
$J$, then $M = P \backslash N$ is a finite set of monomials because
$S/J$ is artinian.  By Lemma \ref{RADICAL1},
$I$ is the ideal of the finite set of points
$\overline{M} \subseteq \pnk$.
Finally, by Proposition
\ref{liftingproperties} the  set $\overline{M}$ is ACM 
because $S/J$ is artinian, and hence, CM.  
\end{proof}

Since characterizing
the Hilbert functions of ACM sets of points in $\pnk$
is equivalent to  characterizing the Hilbert functions of 
$\N^k$-graded artinian quotients of  
$S$, Theorem \ref{acmforpnk} translates one open problem into 
another open problem
because we do not have a theorem like Macaulay's Theorem \cite{Ma} for 
$\N^k$-graded rings if $k > 1$.
However, as shown below, there is a Macaulay-type 
theorem for $\N^k$-graded quotients of  
${\bf k}[x_{1,1},x_{2,1},\ldots,x_{k,1}]$.  As a consequence, we can 
explicitly  describe  all the Hilbert functions of ACM sets of points in 
$\pr^1 \times \cdots \times \pr^1$ ($k$ times) for any positive integer $k$.

So, suppose that $S = {\bf k}[x_1,\ldots,x_k]$ and $\deg x_i = e_i$,
where $e_i$ is the $i^{th}$ standard basis vector 
of $\N^k$.  We prove a stronger result by  characterizing
the Hilbert functions of all quotients of $S$, not only the artinian quotients.

\begin{theorem}	\label{hilbertfornk1}
Let $S = {\bf k}[x_1,\ldots,x_k]$ with
$\deg x_i = e_i$,
and let $H: \N^k \rightarrow \N$ be a numerical function.  Then
there exists a homogeneous ideal $I \subsetneq S$ with Hilbert
function $H_{S/I} = H$ if and only if
\begin{enumerate}
\item[$(i)$] $H(0,\ldots,0) = 1$,
\item[$(ii)$] $H(\bi) = 1$ or $0$ if $\bi >
\underline{0}$, and 
\item[$(iii)$] if $H(\bi) = 0$, then
$H(\bj) = 0$ for all $\bj \geq \bi$.
\end{enumerate}
\end{theorem}

\begin{proof}
If $I \subsetneq S$ is an $\N^k$-homogeneous ideal
such that $H_{S/I} = H$,  then condition
$(i)$ is simply a consequence of the fact that $I \subsetneq S$.  
Statement $(ii)$ follows from the inequality
$0 \leq H_{S/I}(\bi) \leq \dim_{\bf k} S_{\bi} = 1$.
Finally, if
$H_{S/I}(\bi) =0$, then
$x_1^{i_1}\cdots x_k^{i_k} \in I$, or equivalently, $S_{\bi}
\subseteq I$ because $x_1^{i_1}\cdots x_k^{i_k}$ is a 
basis for $S_{\bi}$.  So,
if $\bj \geq \bi$, then
$S_{\bj} \subseteq I$, and hence,
$H_{S/I}(\bj) = 0$, thus proving $(iii)$.

Conversely, suppose that $H$ is a numerical function 
satisfying $(i)-(iii)$.  If $H(\bi) =1$ for
all $\bi \in \N^k$, then the ideal $I = (0) \subseteq
S$ has the property that $H_{S/I} = H$.

So, suppose $H(\bi) \neq 1$ for all $\bi$.
Set $\mathcal{I}:=\left\{ \bi \in \N^k ~\left|~ H(\bi) = 0 
\right\}\right..$ Note that $\mathcal{I} \neq \N^k$ because $\underline{0} 
\not\in \mathcal{I}$. 
Let $I$ be the ideal  $I:=\langle\{x_1^{i_1}\cdots x_k^{i_k} ~|~ \bi 
\in\mathcal{I}\}\rangle$ in $S$.  We claim that $H_{S/I}(\bi) = H(\bi)$
for all $\bi \in \N^k$.  It is immediate that
$H_{S/I}(\underline{0})=H(\underline{0}) =1$.  Moreover,
if $H(\bi) = 0$, then $H_{S/I}(\bi) = 0$ because $x_1^{i_1}\cdots x_k^{i_k} 
\in I_{\bi} \subseteq I$, i.e., $S_{\bi} \subseteq I$.

So, we need to check that  if $H(\bi) =1$,
then $H_{S/I}(\bi) =1$.  Suppose
$H_{S/I}(\bi) = 0$.  This implies that 
$x_1^{i_1}\cdots x_k^{i_k} \in I$.  But because $\bi \not\in
\mathcal{I}$, there is a monomial $x_1^{j_1}\cdots x_k^{j_k} \in I$
with $\bj \in \mathcal{I}$, such that $x_1^{j_1}\cdots x_k^{j_k}$
divides $x_1^{i_1}\cdots x_k^{i_k}$.  But this is equivalent
to the statement that $\bj \leq \bi$.  But
this contradicts hypothesis $(iii)$.  So $H_{S/I}(\bi) = 1$.
\end{proof}

Using Theorem \ref{hilbertfornk1} and the definition
of an $\N^k$-graded artinian quotient we then have:

\begin{corollary}  \label{artiniannk1}
Let $S = {\bf k}[x_1,\ldots,x_k]$ with $\deg x_i = e_i$, and 
let $H:\N^k \rightarrow \N$ be a numerical function.  Then
$H$ is the Hilbert function of an $\N^k$-graded artinian
quotient of $S$ if and only if 
\begin{enumerate}
\item[$(i)$] $H(0,\ldots,0) = 1$,
\item[$(ii)$] $H(\bi) = 1$ or $0$ if $\bi >
(0,\ldots,0)$,
\item[$(iii)$] if $H(\bi) = 0$, then
$H(\bj) = 0$ for all $\bj \geq \bi$, and
\item[$(iv)$] for each $i \in [k]$ there exists
an integer $t_i$ such that $H(t_ie_i) = 0$.
\end{enumerate}
\end{corollary}

\begin{corollary}	\label{p1kcor}
Let $H:\N^k \rightarrow \N$ be a numerical function.  Then $H$
is the Hilbert function of an ACM set of distinct points
in $\pr^1 \times \cdots \times \pr^1$ ($k$ times) if and only
$\Delta H$ satisfies conditions $(i)-(iv)$ of Corollary ~\ref{artiniannk1}.
\end{corollary}

\begin{remark}	\label{p1p1remark}
It follows from the previous corollaries that $H$ is the Hilbert
function of an ACM set of points in $\pr^1 \times \pr^1$ if and
only if 
\begin{enumerate}
\item[$(i)$] $\Delta H(i,j) = 1$ or $0$,
\item[$(ii)$] if $\Delta H(i,j) = 0$, then $\Delta H(i',j') = 0$ for
all $(i',j') \in \N^2$ with $(i',j') > (i,j)$, and
\item[$(iii)$] there exists integers $t$ and $r$ such that
$\Delta H(t,0) = 0$ and $\Delta H(0,r) = 0$.
\end{enumerate}
Giuffrida, Maggioni, and Ragusa proved precisely this result 
in Theorems 4.1 and 4.2 of ~\cite{GuMaRa1}.  We investigate
ACM sets of points in $\popo$ in further detail in the next two sections. 
\end{remark}


\section{ACM Sets of Points in $\popo$}

If $\X$ is an ACM set of points in $\popo$, 
then by Theorem \ref{acmforpnk} the function $\Delta H_{\X}$
is the Hilbert function of a bigraded artinian quotient
of ${\bf k}[x_1,y_1]$.  In \cite{GuMaRa1} it was shown that the
converse of this statement is also true, thereby classifying
the ACM sets of points in $\popo$.  In this section
we revisit this result by giving a new proof of this characterization
that depends only upon numerical information describing $\X$. 

We begin with a brief digression
to introduce some needed combinatorial results.
Recall that a tuple $\lambda = (\lambda_1,\ldots,\lambda_r)$
of positive integers is a {\it partition} of an integer $s$,
denoted $\lambda \vdash s$,
if $\sum \lambda_j = s$ and $\lambda_i \geq \lambda_{i+1}$
for each $i$.
If $\lambda \vdash s$, then the
{\it conjugate} of $\lambda$ is the tuple
$\lambda^* = (\lambda^*_1,\ldots,\lambda_{\lambda_1}^*)$
where $\lambda_i^* := \#\{\lambda_j \in \lambda ~|~
\lambda_j \geq i\}$.  Moreover, $\lambda^*$ is also a partition
of $s$.

To any partition $\lambda = (\lambda_1,\ldots,\lambda_r) \vdash s$
we can associate the following diagram:  on an $r \times \lambda_1$
grid, place $\lambda_1$ points on the first line,
$\lambda_2$ points on the second, and so on.  The resulting diagram is 
called the {\it Ferrers diagram} of $\lambda$.
For example,
suppose $\lambda = (4,4,3,1,1) \vdash 13$.  Then the Ferrers diagram
is
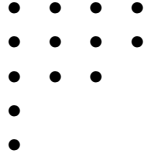
\[
\begin{tabular}{cccc}
$\bullet$ & $\bullet$ & $\bullet$ & $\bullet$ \\
$\bullet$ & $\bullet$ & $\bullet$ & $\bullet$ \\
$\bullet$ & $\bullet$ & $\bullet$ & \\
$\bullet$ & & & \\
$\bullet$ & & &
\end{tabular}\]
The conjugate of $\lambda$ can be read off the Ferrers diagram by
counting the number of dots in each column as opposed to each row.  In
this example $\lambda^* = (5,3,3,2)\vdash 13$.

The following lemma, whose proof is a straightforward combinatorial exercise,
describes some of the relations between a partition and its conjugate.
\begin{lemma}  	\label{conjugateprop}
Let $\alpha = (\alpha_1,\ldots,\alpha_n) \vdash s$, 
and $\beta = (\beta_1,\ldots,\beta_m) \vdash s$.
If $\alpha^* = \beta$, then
\begin{enumerate}
\item[$(i)$] $\alpha_1 = |\beta|$ and $\beta_1 = |\alpha|$.
\item[$(ii)$]  if $\alpha' = (\alpha_2,\ldots,\alpha_n)$ and $\beta'=
(\beta_1 -1,\ldots,\beta_{\alpha_2} - 1)$, then $(\alpha')^* = \beta'$.
\end{enumerate}
\end{lemma}

Let $\X$ denote a set of reduced points in $\popo$, and associate
to $\X$ two tuples $\ax$ and $\bx$ as follows.  Let
$\pi_1(\X) = \{P_1,\ldots,P_t\}$ be the $t$ distinct
first coordinates in $\X$.  Then, for each $P_i \in 
\pi_1(\X)$, let $\alpha_i := |\pi_1^{-1}(P_i)|$, i.e.,
the number of points in $\X$ which have $P_i$ as its
first coordinate.  After relabeling the $\alpha_i$
so that $\alpha_i \geq \alpha_{i+1}$ for $i = 1,\ldots,t-1$,
we set $\ax = (\alpha_1,\ldots,\alpha_t)$.  Analogously,
for each $Q_i \in \pi_2(\X) = \{Q_1,\ldots,Q_r\}$, we let
$\beta_i := |\pi_2^{-1}(Q_i)|$.  After relabeling so that
$\beta_i \geq \beta_{i+1}$ for $i = 1,\ldots,r-1$, we set
$\bx = (\beta_1,\ldots,\beta_r)$.
So, by construction, $\ax$, $\bx \vdash s = |\X|$.
Note that $|\ax| = |\pi_1(\X)|$ and $|\bx| = |\pi_2(\X)|$.

We write the Hilbert function $H_{\X}$ as an infinite matrix
$(m_{ij})$ where $m_{ij}:= H_{\X}(i,j)$.  Proposition \ref{sides} gives
\[
m_{i,0} = \left\{\begin{array}{ll}
i+1 & 0 \leq i \leq t-1 \\
t & i \geq t 
\end{array}
\right.
~~\text{and}~~
m_{0,j} = \left\{\begin{array}{ll}
i+1 & 0 \leq i \leq r-1 \\
r & i \geq r 
\end{array}
\right..
\]
because $\pi_1(\X) = \{P_1,\ldots,P_t\} \subseteq \pr^1$
and $\pi_2(\X) = \{Q_1,\ldots,Q_r\} \subseteq \pr^1$.
This fact, combined with Proposition \ref{eventualgrowth}, 
implies that $H_{\X}$ has the form
\begin{equation}  \label{matrix}
H_{\X} = 
\bmatrix
1& 2 & \cdots  &r-1 & {\bf r} & r & \cdots \\
2&   &         &    &{\bf m_{1,r-1}} & m_{1,r-1} & \cdots \\
\vdots &  & *  &&\vdots & \vdots & \\ 
t-1&  &   && {\bf m_{2,r-1}} & m_{2,r-1} & \cdots \\
{\bf t} & {\bf m_{t-1,1}} & \cdots & {\bf m_{t-1,r-2}}& {\bf s } 
& s & \cdots \\
t & m_{t-1,1} & \cdots  & m_{t-1,r-2} & s  & s \\
\vdots &\vdots & & \vdots& \vdots && \ddots
\endbmatrix
\end{equation}
where the values denoted by $(*)$ need to be calculated.

Set $B_C = (m_{t-1,0},\ldots,m_{t-1,r-1})$ and $B_R = 
(m_{0,r-1},\ldots,m_{t-1,r-1})$.  From our description
of $H_{\X}$, we see that if we know the values in the tuples
$B_C$ and $B_R$, we will know all but a finite
number of values of $H_{\X}$.
As shown below,  
the tuples $B_C$ and $B_R$ can be computed directly from
the tuples $\ax$ and $\bx$ defined above.  If $\lambda$ is 
a tuple, then we shall abuse notation and write $\lambda_i \in \lambda$
to mean that $\lambda_i$ is a coordinate of $\lambda$.

\begin{proposition}[{\cite[Proposition 5.11]{VT2}}] 
\label{computingborderinp1p1}
Let $\X \subseteq \popo$ be a finite set of  points in
$\popo$.  
\begin{enumerate}
\item[$(i)$] if $B_C = (m_{t-1,0},\ldots,m_{t-1,r-1})$ 
where $m_{t-1,j} = H_{\X}(t-1,j)$ for $j = 0,\ldots,r-1$, then
\[m_{t-1,j} = \#\{\alpha_i \in \ax ~|~\alpha_i \geq 1\} +
 \#\{\alpha_i \in \ax ~|~\alpha_i \geq 2\} + \cdots 
+ \#\{\alpha_i \in \ax ~|~\alpha_i \geq j+1\}. \]
\item[$(ii)$] if $B_R = (m_{0,r-1},\ldots,m_{t-1,r-1})$ 
where $m_{j,r-1} = H_{\X}(j,r-1)$ for
$j = 0,\ldots,t-1$, then
\[m_{j,r-1} = \#\{\beta_i \in \bx ~|~\beta_i \geq 1\} +
 \#\{\beta_i \in \bx ~|~\beta_i \geq 2\} + \cdots 
+ \#\{\beta_i \in \bx ~|~\beta_i \geq j+1\}. \]
\end{enumerate}
\end{proposition}

We can rewrite this result more compactly using the 
language of combinatorics  introduced above.
If $p = (p_1,\ldots,p_k)$, then
we write $\Delta p$ to denote the tuple $
\Delta p := (p_1,p_2-p_1,\ldots,p_k - p_{k-1})$.

\begin{corollary}[{\cite[Corollary 5.13]{VT2}}]
Let $\X$ be a finite set of points in $\popo$.  Then 
\begin{enumerate}
\item[$(i)$]  $\Delta B_C = \ax^*$. 
\item[$(ii)$] $\Delta B_R = \bx^*.$
\end{enumerate}
\end{corollary}

\begin{remark}In \cite{VT2} the tuple $B_{\X} = (B_C,B_R)$ was 
called the {\it border} of the Hilbert function. 
\end{remark}

Recall that $\Delta H_{\X}$, the  {\it first difference function}
of $H_{\X}$, is defined by
\[\Delta H_{\X}(i,j) := H_{\X}(i,j) - H_{\X}(i-1,j) - H_{\X}(i,j-1)
+ H_{\X}(i-1,j-1)\]
where $H_{\X}(i,j) = 0$ if $(i,j) \not\geq (0,0)$.
The entries of  $\ax^*$ and $\bx^*$ can then be read
from $\Delta H_{\X}$.

\begin{corollary}
\label{rowcolumnsum}
Let $\X \subseteq \popo$ be a finite set of points, and
set  $c_{i,j} := \Delta H_{\X}(i,j)$.  Then 

\begin{enumerate}
\item[$(i)$] for every $0 \leq j \leq r-1 = |\pixt| -1$
\[\alpha^*_{j+1} = \sum_{h \leq |\pix| -1} c_{h,j}\]
where $\alpha^*_{j+1}$ is the $(j+1)^{th}$ entry of $\ax^*$,
the conjugate of $\ax$.
\item[$(ii)$] for every $0 \leq i \leq t-1 = |\pix| -1$
\[\beta^*_{i+1} = \sum_{h \leq |\pixt| -1} c_{i,h}.\]
where $\beta^*_{i+1}$ is the $(i+1)^{th}$ entry of $\bx^*$.
\end{enumerate}
\end{corollary}

\begin{proof}
Use Proposition \ref{computingborderinp1p1} and the identity
$H_{\X}(i,j) = \sum_{(h,k) \leq (i,j)} c_{h,k}$
to compute $\alpha^*_{j+1}$: 
\begin{eqnarray*} 
\alpha^*_{j+1} & = & H_{\X}(t-1,j) - H_{\X}(t-1,j-1) \\
& = & \sum_{(h,k)\leq (t-1,j)} c_{h,k} - 
\sum_{(h,k) \leq (t-1, j-1)} c_{h,k} = 
\sum_{h \leq t-1=|\pix|-1} c_{h,j}.
\end{eqnarray*}
The proof for the second statement is the same.
\end{proof}

\begin{lemma}	\label{projectionLemma}
Let $\X \subseteq \popo$ be a finite set of points, and suppose that 
$\ax^* = \bx$.  Let $P$ be 
a point of $\pix$ such that $|\pi_1^{-1}(P)| = \alpha_1$.  Set 
$\X_P := \pi_1^{-1}(P)$.  Then $\pi_2(\X_P) = \pi_2(\X).$
\end{lemma}

\begin{proof}
Since $\X_P \subseteq \X$, we have
$\pi_2(\X_P) \subseteq \pixt$.   Now, by our choice of $P$,
$|\pi_2(\X_P)| = \alpha_1$.  But since $|\pixt| = |\bx|$ and
$\ax^* = \bx$, from Lemma ~\ref{conjugateprop}
we have $|\pixt| = |\bx| = \alpha_1 = |\pi_2(\X_{P})|,$ and hence
$\pi_2(\X_P) = \pi_2(\X).$
\end{proof}

\begin{proposition}  \label{ciresolution}
Suppose that $\X$ is a set of $s= tr$ points in $\popo$ such that 
$\ax=(r,\ldots,r)$ ($t$ times) and $\bx = ({t,\ldots,t})$ ($r$ times).
Then $\X$ is a complete intersection, and 
the graded minimal free resolution of $\Ix$ is given by
\[0 \longrightarrow R(-t,-r) \longrightarrow 
R(-t,0)\oplus R(0,-r) \longrightarrow \Ix \longrightarrow 0.\]
\end{proposition}

\begin{proof}
Because $|\ax|=t$ and $|\bx| = r$, $\pi_1(\X) = 
\{P_1,\ldots,P_t\}$ and $\pi_2(\X) =\{Q_1,\ldots,Q_r\}$ where
$P_i,Q_j \in \pr^1$.  Since $|\X| = tr$, 
$\X =
\left\{P_i \times Q_j \left|~ 1\leq i \leq t, 1\leq j\leq r \right\}
\right..$
Hence, if $I_{P_i\times Q_j} = (L_{P_i},L_{Q_j})$ is the
 ideal associated to the point $P_i \times Q_j$, then the
defining ideal of $\X$ is 
\[
\Ix = \bigcap_{i,j}
\left(L_{P_i},L_{Q_j}\right)
=
\left(L_{P_1}L_{P_2}\cdots L_{P_t},L_{Q_1}L_{Q_2}\cdots L_{Q_r}\right).
\]
Since $\deg L_{P_1}L_{P_2}\cdots L_{P_t} = (t,0)$ and
$\deg L_{Q_1}L_{Q_2}\cdots L_{Q_r} = (0,r)$, the two
generators  form a regular sequence on $R$, and hence,
$\X$ is a complete intersection.  
The graded 
minimal
free resolution is then given by the {\it Koszul resolution},
taking into consideration that $\Ix$ is bigraded.
\end{proof}

We now come to the main result of this section.

\begin{theorem}	\label{acmp1p1}
Let $\X$ be a finite set of  points in $\popo$ 
with Hilbert function $H_{\X}$.   Then
the following are equivalent:
\begin{enumerate}
\item[$(i)$] $\X$ is ACM.
\item[$(ii)$] $\Delta H_{\X}$ is the Hilbert function
of an $\N^2$-graded artinian quotient of ${\bf k}[x_1,y_1]$.
\item[$(iii)$] $\ax^* = \bx$.
\end{enumerate}
\end{theorem}

\begin{proof}
The implication $(i) \Rightarrow (ii)$ is Corollary \ref{ahilbert2}.
So, suppose that $(ii)$ holds.  Because $\Delta H_{\X}$
is the Hilbert function of an $\N^2$-graded artinian quotient
of ${\bf k}[x_1,y_1]$, Corollary \ref{p1kcor}, 
Remark \ref{p1p1remark}, 
and the matrix (\ref{matrix}) give
\begin{center}
\begin{picture}(150,100)(-20,30)
\put(-60,75){$\Delta H_{\X} =$}
\put(10,25){\line(0,1){115}}
\put(10,144){${\scriptstyle 0}$}
\put(0,130){${\scriptstyle 0}$}
\put(86,144){${\scriptstyle r-1}$}
\put(-17,40){${\scriptstyle t-1}$}
\put(10,140){\line(1,0){110}}
\put(10,40){\line(1,0){20}}
\put(30,40){\line(0,1){40}}
\put(30,80){\line(1,0){40}}
\put(70,80){\line(0,1){40}}
\put(70,120){\line(1,0){25}}
\put(95,120){\line(0,1){20}}
\put(35,105){{\bf 1}}
\put(80,50){{\bf 0}}
\end{picture}
\end{center}
where $t = |\pi_1(\X)|$ and $r = |\pi_2(\X)|$.
We have written $\Delta H_{\X}$ as an infinite matrix whose indexing starts
from zero rather than one.  

By Corollary \ref{rowcolumnsum} the number of $1$'s in the $(i-1)^{th}$ row
of $\Delta H_{\X}$ for each integer $1 \leq i \leq t$
is simply the $i^{th}$ coordinate of $\bx^*$.  Similarly, the number of ones 
in the $(j-1)^{th}$ column of $\Delta H_{\X}$ for each integer
$1 \leq j \leq r$ is the $j^{th}$ coordinate of $\ax^*$. Now $\Delta H_{\X}$
can be identified with the Ferrers diagram 
 of $\bx^*$ by associating to each $1$ in $\Delta H_{\X}$ 
a dot in the Ferrers diagram in the natural way, i.e.,

\begin{picture}(150,150)(-40,10)
\put(10,25){\line(0,1){115}}
\put(10,144){${\scriptstyle 0}$}
\put(0,130){${\scriptstyle 0}$}
\put(86,144){${\scriptstyle r-1}$}
\put(-17,40){${\scriptstyle t-1}$}
\put(10,140){\line(1,0){110}}
\put(10,40){\line(1,0){20}}
\put(30,40){\line(0,1){40}}
\put(30,80){\line(1,0){40}}
\put(70,80){\line(0,1){40}}
\put(70,120){\line(1,0){25}}
\put(95,120){\line(0,1){20}}
\put(35,105){{\bf 1}}
\put(170,50){$\bullet$}
\put(170,70){$\bullet$}
\put(170,90){$\bullet$}
\put(170,130){$\bullet$}
\put(170,110){$\bullet$}
\put(190,90){$\bullet$}
\put(190,110){$\bullet$}
\put(190,130){$\bullet$}
\put(210,90){$\bullet$}
\put(210,110){$\bullet$}
\put(210,130){$\bullet$}
\put(230,130){$\bullet$}
\put(130,90){$\longleftrightarrow$}
\end{picture}

\noindent
By using the Ferrers diagram and Corollary \ref{rowcolumnsum}, it is now 
straightforward to calculate that the conjugate of $\bx^*$ is 
$(\bx^*)^* = \bx = \ax^*$, and so $(iii)$ holds.

To demonstrate that $(iii)$ implies $(i)$, we proceed 
by induction on the tuple $(|\pix|,|\X|)$.  For any positive integer $s$,
if $(|\pix|,|\X|) = (1,s)$, then
$\ax = (s)$ and $\bx = ({1,\ldots,1})$ ($s$ times), and so
$\ax^* = \bx$. Then by Proposition \ref{ciresolution} $\X$ is 
a complete intersection, and hence, ACM.

So, suppose that $(|\pi(\X)),|\X|) = (t,s)$ and that 
result holds true for all $\Y \subseteq \popo$
with $\alpha_{\Y}^* = \beta_{\Y}$ and $(t,s) >_{lex} (|\pi_1(\Y)|,|\Y|)$
 where $>_{lex}$ is the lexicographical ordering 
on $\N^2$.

Suppose that $P_1$ (after a possible relabeling) is the
element of $\pix$ such that $|\pi_1^{-1}(P_1)| = \alpha_1$.  Let
$L_{P_1}$ be the form of degree $(1,0)$  that vanishes at $P_1$.  By
abusing notation, we also let $L_{P_1}$ denote the
$(1,0)$-line in $\popo$ defined by $\lp$. 

Set $\xp := \X \cap \lp = \pi_1^{-1}(P_1)$ and $\Z := \X \backslash \xp$.  
It follows that $\alpha_{\Z} = (\alpha_2,\ldots,\alpha_t)$
and $\beta_{\Z} =(\beta_1 -1,\ldots,\beta_{\alpha_2} -1)$.   
Now $(t,s) >_{lex} (|\pi_1(\Z)|,|\Z|)$, and  moreover, 
 $\alpha_{\Z}^* = \beta_{\Z}$ by Lemma \ref{conjugateprop}.
Thus, by the induction hypothesis, $\Z$ is ACM.

Suppose that $\pixt = \{Q_1,\ldots,Q_r\}$.  Let $L_{Q_i}$ be
the degree $(0,1)$ form that vanishes at $Q_i \in \pixt$ and
set $F: = L_{Q_1}L_{Q_2}\cdots L_{Q_r}$.  
Because $\ax^* = \bx$, from Lemma ~\ref{projectionLemma} we have
$\pi_2(\xp) = \pixt$.  So, $\X_{P_1} = 
\left\{P_1\times Q_1,\ldots,P_1\times Q_r\right\}$, and hence
$
\Ixp = \bigcap^r_{i=1} (L_{P_1},L_{Q_i}) = (L_{P_1},F).
$
Furthermore, if $P \times Q \in \Z$, then $Q \in \pi_2(\Z) \subseteq \pixt$,
and thus $F(P \times Q) = 0$.  Therefore $F \in I_{\Z}$.
Because $F$ is in $\Iz$ and is also a generator of $\Ixp$, we 
can show:

{\it Claim. }  Let $I = L_{P_1}\cdot\Iz + (F).$  Then $I = \Ix$.

{\it Proof of the Claim. }  Since $\Ix = I_{\Z \cup \xp} = 
\Iz \cap \Ixp,$ we will show $\Iz \cap \Ixp =
 L_{P_1}\cdot\Iz + (F)$.
So, suppose that $G = L_{P_1}H_1 + H_2F \in  L_{P_1}\cdot\Iz + (F)$
with $H_1 \in \Iz$ and $H_2 \in R$.  Because
$L_{P_1}$ and $F$ are in $\Ixp$,  we have $G \in \Ixp$.
Since $H_1,F \in \Iz$, $G \in \Iz$.  Thus $G \in I_{\Z} \cap \Ixp$.

Conversely, let $G \in \Iz \cap \Ixp.$  Since $G \in
\Ixp$,  $G = \lp H_1 + F H_2$. We need to show that $H_1 \in \Iz$.
Because
$G,F \in \Iz$, we also have $\lp H_1 \in \Iz$.  But for every
$P \times Q \in \Z$, $P \neq P_1$, and thus $\lp(P\times Q) \neq 0$.
Hence $\lp H_1 \in \Iz$ if and only if $H_1(P \times Q) = 0$
for every $P \times Q \in \Z$.
\hfill$\Box$

We note that 
$\X \subseteq \popo$ is ACM if and only if
the variety $\widetilde{\X} \subseteq \pr^3$ defined by $\Ix$,
considered as a homogeneous ideal of $R = {\bf k}[x_0,x_1,y_0,y_1]$, is ACM.
As a variety of $\pr^3$,  $\widetilde{\X}$ is a curve since
K-$\dim R/\Ix  = 2$.  Let $\widetilde{\Z}$
denote the curve of $\pr^3$ defined by $\Iz$, considered
also as a homogeneous ideal of $R$.  The
claim implies that the curve $\widetilde{\X}$ is a basic double
link of 
$\widetilde{\Z}$.  Since the Cohen-Macaulay property is preserved under 
linkage (see \cite[Theorem 3.2.3]{M} and following remark),
$\widetilde{\X}$ is an ACM curve of $\pr^3$, or equivalently, $\X$ is
an ACM set of points in $\popo$.
\end{proof}

\begin{remark}
Giuffrida, {\it et al.} \cite[Theorem 4.1]{GuMaRa1}
demonstrated the equivalence of statements $(i)$ and $(ii)$
of Theorem \ref{acmp1p1} via different means. Our
contribution is to show that the ACM sets of points are also
characterized  by   $\ax$ and $\bx$.  This
result has been extended in \cite{GVT} to characterize ACM fat point schemes
in $\popo$ 
\end{remark}

\begin{corollary} \label{removeptfromacm}
Let $\X$ be a set of points in $\popo$ with 
$\ax = (\alpha_1,\ldots,\alpha_t)$, and $\pi_1(\X) = \{P_1,\ldots,
P_t\}$.  Suppose (after a possible relabeling) that 
$|\pi_1^{-1}(P_i)| = \alpha_i$.  
Set
\[
\X_i := \X \backslash \left\{ \pi_1^{-1}(P_1) \cup \cdots \cup 
\pi_1^{-1}(P_i) \right\} \hspace{.5cm}\mbox{for $0 \leq i \leq t-1$},
\]
where $\X_0 := \X$.  If $\X$ is ACM, then, for each integer 
$0 \leq i \leq t-1$,   $\X_i$ is ACM 
with $\alpha_{\X_i} = (\alpha_{i+1},
\alpha_{i+2},\ldots,\alpha_t)$.
\end{corollary}

\begin{proof}
It is sufficient to show that for each $i = 0,\ldots,t-2$,
if $\X_i$ is ACM, then $\X_{i+1}$ is ACM.  
Since $\X_{i+1} = \X_i \backslash \{\pi_1^{-1}(P_{i+1})\},$
$\X_{i+1}$ is constructed from $\X_i$ by removing the $\alpha_{i+1}$
points of $\X_i$ which have $P_{i+1}$ as its first coordinate.  
The tuple $\beta_{\X_{i+1}}$ is constructed from $\beta_{\X_i}$
by subtracting $1$ from $\alpha_{i+1}$ coordinates in $\beta_{\X_i}$.
But  because $\alpha_{\X_i}^* = \beta_{\X_i}$, we have $|\beta_{\X_i}| 
= \alpha_{i+1}$, and 
thus $\beta_{\X_{i+1}} = (\beta_1-1 ,\ldots,\beta_{\alpha_{i+1}}-1)
= (\beta_1-1,\ldots,\beta_{\alpha_{i+2}}-1)$.  But by
Lemma \ref{conjugateprop}, $\alpha_{\X_{i+1}}^* = \beta_{\X_{i+1}}$, and 
hence, $\X_{i+1}$ is ACM by Theorem \ref{acmp1p1}.
\end{proof}


\section{The Betti numbers of ACM sets of points in $\popo$}
 
If $\X$ is a set of $s$ points in $\pr^1$, then it is well known
that the graded Betti numbers,
and consequently, the Hilbert function of $\X$, can be determined
solely from $|\X| = s$.  If we restrict to ACM sets of points in $\popo$,
we can extend this result to show that the graded Betti numbers
in the minimal free resolution of $\X$ (and hence, $H_{\X}$) 
can be computed directly from the tuples $\ax$ and $\bx$ introduced
in the previous section which
numerically describe  $\X$.

We require
the following notation to describe the minimal free resolution. 
 Suppose that $\X \subseteq \popo$
is a set of points with $\ax = (\alpha_1,\ldots,\alpha_t)$.  Define
the following sets:
\begin{eqnarray*}
C_{\X} & := &\left\{(t,0),(0,\alpha_1)\right\} \cup
\left\{(i-1,\alpha_i) ~|~ \alpha_i - \alpha_{i-1} < 0\right\}, \\
V_{\X}& := & \left\{ (t,\alpha_t) \right\} \cup 
\left\{ (i-1,\alpha_{i-1}) ~|~ \alpha_i-\alpha_{i-1} < 0 \right\}.
\end{eqnarray*}
We take $\alpha_{-1} = 0$.  With this notation, we have

\begin{theorem}  \label{acmresolution}
Suppose that $\X$ is an ACM set of points in $\popo$.  
Let $C_{\X}$ and $V_{\X}$
be constructed from $\ax$ as above.  Then
the graded  minimal free resolution of $\Ix$ is given by
\[
0 \longrightarrow \bigoplus_{(v_1,v_2) \in V_{\X}} R(-v_1,-v_2) 
\longrightarrow \bigoplus_{(c_1,c_2) \in C_{\X}} R(-c_1,-c_2) 
\longrightarrow \Ix \longrightarrow 0.
\]
\end{theorem}

\begin{proof}
We proceed by induction on the tuple $(|\pi_1(\X)|,|\X|)$.  
If $s$ is any integer, and $(|\pi_1(\X),|\X|) = (1,s)$, then $\ax = (s)$ and 
$\bx = (1,\ldots,1)$ ($s$ times).  The conclusion follows from 
Proposition \ref{ciresolution} since $C_{\X} = \{(1,0),(0,s)\}$ and 
$V_{\X} = \{(1,s)\}$.

So, suppose $(|\pi_1(\X),|\X|) = (t,s)$ with $t > 1$.
Then $\ax = (\underbrace{\alpha_1,\ldots,\alpha_1}_l,
\alpha_{l+1},\ldots,\alpha_t)$, i.e., $\alpha_{l+1} < \alpha_1$, but 
$\alpha_l = \alpha_1$.
If $l=t$, then $\X$ is a complete intersection and the resolution
is given by Proposition \ref{ciresolution}.  The conclusion
now follows because $C_{\X} = \{(l,0),(0,\alpha_1)\}$ and $V_{\X} = 
\{(l,\alpha_1)\}$.

If  $l < t$,  let $P_1, \ldots,P_l$ be the $l$ points
of $\pi_1(\X)$ that have $|\pi_1^{-1}(P_i)| = \alpha_1$.  Set
$\Y = \pi_1^{-1}(P_1) \cup \cdots \cup \pi_1^{-1}(P_l)$.  Because
$\X$ is ACM, $\alpha_1 = |\bx|$, and hence,
$\Y = \left\{ P_i \times Q_j \left|~ 1\leq i \leq l,~ Q_j \in \pi_2(\X)
\right\}\right.$.
So, $\alpha_{\Y} = ({\alpha_1,\ldots,\alpha_1})$ and
$\beta_{\Y} = ({l,\ldots,\l})$,
and thus $\Y$ is a complete intersection.  In fact,
$I_{\Y} = (L_{P_1}\cdots L_{P_l},L_{Q_1}\cdots L_{Q_{\alpha_1}})$
where $L_{P_i}$ is the form of degree $(1,0)$ that vanishes at all 
the points of
$\popo$ which have $P_i$ as their first coordinate, and $L_{Q_i}$
is the form of degree $(0,1)$ that vanishes at all points $P\times Q \in \popo$
such that $Q = Q_i$.  

Let $F: =  L_{P_1}\cdots L_{P_l}$ and $G := L_{Q_1}\cdots L_{Q_{\alpha_1}}$.
By Proposition \ref{ciresolution} the resolution of $I_{\Y}$ is 
\[0 \longrightarrow R(-l,-\alpha_{1}) 
\stackrel{\phi_2}{\longrightarrow} 
R(-l,0)\oplus R(0,-\alpha_1) 
\stackrel{\phi_1}{\longrightarrow} I_{\Y} \longrightarrow 0
\]
where $\phi_1 = [F ~~ G]$
and $\phi_2 = 
\bmatrix
G \\
-F
\endbmatrix$.
Let $\Z := \X \backslash \Y$.  Since $\pi_2(\Z) \subseteq \pi_2(\X)$,
it follows that $G = L_{Q_1}\cdots L_{Q_{\alpha_1}} \in I_{\Z}$.  
Hence, $\operatorname{im} \phi_2 \subseteq I_{\Z}(-l,0) 
\oplus R(0,-\alpha_1)$.  

{\it Claim. }  $\Ix = F\cdot I_{\Z} + (G)$

{\it Proof of the Claim. }  By construction, $\X = \Z \cup \Y$.  
Hence, we want to show that
$I_{\Z} \cap I_{\Y} = F\cdot I_{\Z} + (G)$.  The proof
now follows as in the proof of the Claim in Theorem 
\ref{acmp1p1}.\hfill $\Box$

From the above resolution for $I_{\Y}$, the claim, and
the fact that  $\operatorname{im} \phi_2 \subseteq I_{\Z}(-l,0) 
\oplus R(0,-\alpha_1)$, we have the following 
 short exact sequence of $R$-modules
\[0 \longrightarrow R(-l,-\alpha_1) 
\stackrel{\phi_2}{\longrightarrow} 
I_{\Z}(-l,0)\oplus R(0,-\alpha_1) 
\stackrel{\phi_1}{\longrightarrow} \Ix = F\cdot I_{\Z} + (G)
\longrightarrow 0
\]
where $\phi_1$ and $\phi_2$ are as above.  

By Corollary \ref{removeptfromacm} the set  $\Z$
is ACM  with $\alpha_{\Z} =
(\alpha_{l+1},\ldots,\alpha_t)$.  Therefore, the induction
hypothesis holds for $\Z$.  With the above short exact 
sequence, we can use the {\it mapping cone construction} 
(see Section 1.5 of \cite{W}) 
to construct a resolution for $\Ix$.  In particular,
we get
\[0 \longrightarrow
\left[\bigoplus_{(v_1,v_2) \in V_{\Z}} R(-(v_1+l),-v_2)\right] 
\oplus R(-l,-\alpha_1)\longrightarrow \]
\[
\left[\bigoplus_{(c_1,c_2) \in C_{\Z}} R(-(c_1+l),-c_2)\right]
\oplus R(0,\alpha_1)
\longrightarrow \Ix
\longrightarrow 0.
\]
Since the resolution has length 2, and because $\X$ is ACM, the resolution
of $\Ix$ cannot be made shorter by the Auslander-Buchsbaum formula
(cf. \cite[Theorem 4.4.15]{W}).

 To show that this resolution
is minimal, it is enough to show that no tuple in
the set $ \left\{(c_1+l,c_2) \left|~ (c_1,c_2) \in C_{\Z}\right\}\right.
\cup \{(0,\alpha_1)\}$ is in the set 
$\left\{ (v_1+l,v_2)\left|~ (v_1,v_2) 
\in V_{\Z}\right\}\right.\cup \{(l,\alpha_1)\}$.  By the induction
hypothesis, we can assume that no $(c_1,c_2) \in C_{\Z}$
is in $V_{\Z}$, and hence, if $(c_1 +l,c_2) \in 
\left\{(c_1+l,c_2) \left|~ (c_1,c_2) \in C_{\Z}\right\}\right.$, then
$(c_1 +l,c_2)$ is not in 
$\left\{ (v_1+l,v_2)\left|~ (v_1,v_2) \in V_{\Z}\right\}\right.$
If $(c_1+l,c_2) = (l,\alpha_1)$ for some $(c_1,c_2)\in C_{\Z}$,
then this implies that $(0,\alpha_1)$.  But this contradictions
the induction hypothesis.  Similarly, if $(0,\alpha_1) \in
\left\{ (v_1+l,v_2)\left|~ (v_1,v_2) 
\in V_{\Z}\right\}\right.$, this implies $(-l,\alpha_1) \in V_{\Z}$,
which is again a contradiction of the induction hypothesis.  So the
resolution  above is  minimal.

To complete the proof we only need to verify that
\begin{enumerate}
\item[$(i)$]
$C_{\X} = \left\{(c_1+l,c_2) \left|~ (c_1,c_2) \in C_{\Z}\right\}\right.
\cup \{(0,\alpha_1)\}$.
\item[$(ii)$]
$V_{\X} = \left\{ (v_1+l,v_2)\left|~ (v_1,v_2) 
\in V_{\Z}\right\}\right.\cup \{(l,\alpha_1)\}$.
\end{enumerate}
Because the verification of these statements is tedious, but
elementary, we omit the details.
\end{proof}

\begin{remark}
It was shown in \cite[Theorem 4.1]{GuMaRa1}   that the graded Betti numbers
for an ACM set of points $\X \subseteq \popo$ could be determined
via the first difference function $\Delta H_{\X}$, i.e.,
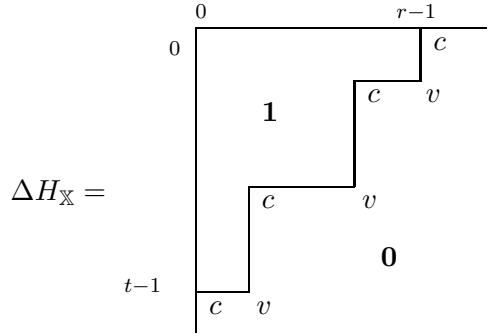
\begin{center}
\begin{picture}(150,150)(-40,0)
\put(-60,75){$\Delta H_{\X} =$}
\put(10,25){\line(0,1){115}}
\put(10,144){${\scriptstyle 0}$}
\put(0,130){${\scriptstyle 0}$}
\put(86,144){${\scriptstyle r-1}$}
\put(-17,40){${\scriptstyle t-1}$}
\put(10,140){\line(1,0){110}}
\put(10,40){\line(1,0){20}}
\put(30,40){\line(0,1){40}}
\put(30,80){\line(1,0){40}}
\put(70,80){\line(0,1){40}}
\put(70,120){\line(1,0){25}}
\put(95,120){\line(0,1){20}}
\put(35,105){{\bf 1}}
\put(80,50){{\bf 0}}
\put(15,32){$c$}
\put(35,72){$c$}
\put(75,112){$c$}
\put(100,132){$c$}
\put(33,32){$v$}
\put(73,72){$v$}
\put(97,112){$v$}
\end{picture}
\end{center}
An element of $C_{\X}$, which \cite{GuMaRa1} 
called a {\it corner} of $\Delta H_{\X}$,
corresponds to a tuple $(i,j)$ that is either 
$(t,0), (0,\alpha_1)=(0,r)$, or has the property that
$\Delta H_{\X}(i,j)=0$, but $\Delta H_{\X}(i-1,j) = \Delta H_{\X}(i,j-1)
=1$.  We have labeled the corners of $\Delta H_{\X}$ with a $c$ in the
above diagram.  An element of $V_{\X}$ is a {\it vertex}.
A tuple $(i,j)$ is called a vertex if $\Delta H_{\X}(i,j) =
\Delta H_{\X}(i-1,j) = \Delta H_{\X}(i,j-1) = 0$, but
$\Delta H_{\X}(i-1,j-1) = 1$.  We have labeled the vertices of $\Delta H_{\X}$
with a $v$ in the above diagram.  Besides giving a new proof
for the resolution of an ACM set of points in $\popo$, we have  shown
that the graded Betti numbers can be computed directly from the tuple $\ax$.
\end{remark}

Using the resolution as given in Theorem \ref{acmresolution} we
can compute $H_{\X}$ directly from $\ax$ for any ACM
set of points $\X \subseteq \popo$.  Formally

\begin{corollary}  \label{computehilbertofacm}
Let $\X$ be an ACM set of points in $\popo$ with $\ax = (\alpha_1,\ldots
\alpha_t)$.  Then
\begin{eqnarray*}
H_{\X} & = &  \bmatrix
1& 2 & \cdots  & \alpha_{1}-1 & \alpha_1 & \alpha_1 & \cdots \\
1& 2 &\cdots  & \alpha_{1}-1 & \alpha_1 & \alpha_1 & \cdots \\
\vdots & \vdots && \vdots  & \vdots &\vdots&\ddots
\endbmatrix
+ 
\bmatrix
0& 0 & \cdots  & 0 & 0 & 0 & \cdots\\
1& 2 & \cdots  & \alpha_{2}-1 & \alpha_2 & \alpha_2 & \cdots \\
1& 2  &\cdots  & \alpha_{2}-1 & \alpha_2 & \alpha_2 & \cdots \\
\vdots & \vdots && \vdots  & \vdots &\vdots&\ddots
\endbmatrix
+ \\
& &  
\bmatrix
0& 0 & \cdots  & 0 & 0 & 0 &\cdots\\
0& 0 & \cdots  & 0 & 0 & 0 &\cdots\\
1& 2 & \cdots  & \alpha_{3}-1 & \alpha_3 & \alpha_3 & \cdots \\
1& 2  &\cdots  & \alpha_{3}-1 & \alpha_3 & \alpha_3 & \cdots \\
\vdots & \vdots & &\vdots  & \vdots &\vdots& \ddots
\endbmatrix
+ \cdots + 
\bmatrix
0& 0 & \cdots  & 0 & 0 & 0 &\cdots\\
\vdots& \vdots& & \vdots & \vdots & \vdots \\
0& 0 & \cdots  & 0 & 0 & 0 &\cdots\\
1& 2 & \cdots  & \alpha_{t}-1 & \alpha_t & \alpha_t & \cdots \\
1& 2  &\cdots  & \alpha_{t}-1 & \alpha_t & \alpha_t & \cdots \\
\vdots & \vdots & &\vdots  & \vdots &\vdots&\ddots
\endbmatrix.
\end{eqnarray*}
\end{corollary}


\section*{acknowledgments}
Some of the results in this paper were part of \cite{VT}.  
The computer program CoCoA \cite{C} was invaluable
in providing and checking examples in the early stages.  
I would like to thank A. Ragusa, L. Roberts, D. Wehlau, and the audience of 
the Curve Seminar for their comments and suggestions.
I would also like to thank the  Universit\`a di Genova,
where part of this work was completed.
I would especially like to thank my supervisor
Tony Geramita for introducing me to this problem, and for his
encouragement and help.


\end{document}